\documentclass[11pt]{amsart}
\usepackage[margin=3cm]{geometry}
\title[Random walks generated by Ewens dist. on the symmetric group]{Random walks generated by the Ewens distribution on the symmetric group} 
\author{Alperen  \"{O}zdemir}
\address{Department of Mathematics, Georgia Institute of Technology} 
\email{aozdemir6@gatech.edu} 

\usepackage{amssymb, amsmath, amsthm, enumerate, mathtools}
\usepackage{float, graphicx, mathtools, pgfplots, tikz, ytableau, xcolor, genyoungtabtikz}
\usetikzlibrary{intersections}

\usepackage[colorlinks=true, urlcolor=black,linkcolor=blue, citecolor=blue]{hyperref}

\newtheorem{theorem}{Theorem}[section]
\newtheorem{lemma}{Lemma}[section]
\newtheorem{corollary}{Corollary}[section]
\newtheorem{definition}{Definition}[section]

\newtheorem{conjecture}{Conjecture}[section]

\newcommand{\stirling}[2]{\genfrac{[}{]}{0pt}{}{#1}{#2}}
\newcommand{\bbox}{\hfill $\Box$}
\newcommand{\pf}{\noindent {\it Proof:} }

\DeclarePairedDelimiter{\floor}{\lfloor}{\rfloor}

\keywords{Ewens distribution, mixing time, random walks on groups, YJM elements}
\subjclass[2020]{60C05, 05E10}

\begin{document}

\begin{abstract}
This paper studies Markov chains on the symmetric group $S_n$ where the transition probabilities are given by the Ewens distribution with parameter $\theta>1$. The eigenvalues are identified to be proportional to the content polynomials of partitions. We show that the mixing time is bounded above by a constant depending only on the parameter if $\theta$ is fixed. However, if it agrees with the number of permuted elements ($\theta=n$), the sequence of chains has a total variation cutoff at $\frac{\log n}{\log 2}.$ 
\end{abstract}

\maketitle

\section{Introduction}
The Ewens distribution on the symmetric group is defined by probabilities that are exponentially proportional to the number of cycles of permutations, with a base parameter $\theta > 0. $ It is originated in population genetics \cite{E}, then has gained a broader reach. It is initially defined on partitions of integers,  which is related to the symmetric group by cycle decompositions. It arises from well-known cycle generating processes, Feller Coupling and Chinese restaurant process. (See \cite{ABT}.) Crane \cite{Cr} has an extensive survey on Ewens sampling formula. 

The Ewens distribution can also be considered a distortion of the uniform distribution on the symmetric group. Hanlon \cite{H} studied a Metropolis Markov chain driven by random transpositions with stationary distribution as the Ewens measure. The spectral analysis of this chain is further carried out in \cite{DH}, and the mixing time is studied. More recently, Jiang \cite{Ji} showed a cutoff result for the random transposition Metropolis chain. 

The chains we study in this paper have transition probabilities given by the Ewens distribution where the parameter $\theta$ is allowed to be a function $n$ for the symmetric group $S_n$. At each stage, a permutation is randomly chosen and multiplied by the permutation of the current state. We are interested in finding the mixing time of the chain depending on the base parameter of the distribution. A more precise description is as follows:  

 Let $\sigma$ and $\tau$ be two permutations in $S_n$. Define $\gamma(\sigma)$ to be the number of cycles in $\sigma.$ We choose $\theta > 1,$ possibly a function of $n$, and raise it to the number of cycles of permutations to obtain the probabilities. So, if the chain is at permutation $\sigma$ at a certain stage, the probability of moving  to $\tau$ at the next stage is
\begin{equation} \label{MC}
 P_{\theta}(\tau\sigma^{-1})= \frac{\theta^{\gamma( \tau\sigma^{-1})} }{\theta^{(n)}},
\end{equation}
where $\frac{1}{\theta^{(n)}}$ is the normalizing constant. Explicitly, 
\begin{equation*}
\theta^{(n)}=\sum_{i=0}^n \stirling{n}{i} \theta^i=\theta (\theta+1)\cdots (\theta+n-1),
\end{equation*}
where $\stirling{n}{i}$ is a Stirling number of the first kind, which counts the number of permutations of cycle length $i$ in $S_n$ (See Section 1.4 of \cite{A}). 

Let us briefly mention some properties of our Markov chain. The chain is irreducible and aperiodic, since the transition probability is non-zero for any permutation. It is time-homogenous, simply because we have the same transition probabilities regardless of the stage. It is also reversible, following from the fact that any permutation and its inverse, which have the same cycle type, have the same probability of being selected at any stage. From the very same fact, the transition matrix is also symmetric, which implies it is doubly stochastic. Therefore, the stationary distribution is the uniform distribution over $S_n.$

Next we define the distance concept used throughout the paper. Let the probability assigned to  $\sigma \in S_n$  be $P_{\theta}^{t}(\sigma)$ after running the chain for $t$ steps. The total variation distance between $P_{\theta}^{t}$ and the uniform distribution $\pi$ is
\begin{equation} \label{TV}
d_{\theta}(t):=\| P_{\theta}^{t} - \pi \|_{TV} = \frac{1}{2} \sum_{\sigma \in S_n} |P_{\theta}^{t}(\sigma) - \pi(\sigma)| = \max_{S \subseteq S_n}|P_{\theta}^{t}(S) - \pi(S)|.
\end{equation}
Our goal is to investigate the convergence time with respect to the total variation distance. It is a common phenomenon in this setting that the convergence to the stationary distribution shows a \textit{cutoff}, a sharp transition. The following definition of cutoff is found in \cite{LPW}. A sequence of Markov chains exhibits a cutoff at the \textit{mixing time} $\{t_n\}$ with a \textit{window} of size $\{w_n\}$ if
\begin{itemize}
\item[(i)] $\lim_{n\rightarrow \infty} \frac{w_n}{t_n} = 0,$
\item[(ii)] $\lim_{\alpha \rightarrow - \infty} \liminf_{n\rightarrow \infty } d_{TV}(t_n + \alpha w_n) = 1,$
\item[(iii)] $\lim_{\alpha \rightarrow  \infty} \limsup_{n\rightarrow \infty } d_{TV}(t_n + \alpha w_n) = 0.$
\end{itemize}
Also, see \cite{Dc} for an insigthful treatment of the concept and \cite{S} for another definition along with a long list of examples. Our result reads:

 \label{main}\begin{theorem}  Let $P_{\theta}$ be the Markov chain on $S_n$ defined above in \eqref{MC}. 
\begin{itemize}
\item[(i)] If $\alpha, \theta>1$ and $t \geq \alpha \frac{\theta^2+1}{2},$ 
then $d_{\theta}(t) < C(\theta)e^{-\alpha}$ for some constant $C(\theta).$ 

\vspace*{1mm}

\item[(ii)] If $\theta = n$, then the chains have a total variation cutoff at $t=\frac{\log n}{\log 2}.$
\end{itemize}
\end{theorem}

For the rest of the paper, we first present the representation theory tools in Section \ref{chp1}. Then we identify the eigenvalues of the Markov chain using Young-Jucys-Murphy elements. They turn out to be functions of contents of Young diagrams. In the following section, we prove the first part of Theorem \ref{main} by using Plancherel growth process. We also provide a simpler proof for only the integer values of $\theta.$ In Section \ref{pf2}, we prove the second part of the theorem.

\section{Representation theory techniques} \label{chp1}

We summarize the techniques that are used in the following sections. We start with the method developed in \cite{DS} to study the random transposition walk. For the connection to the group representation theory, first consider the \textit{Fourier transform} of the measure $\mu$ defined on a finite group $G$, evaluated at any representation $\lambda$ of $G$, 

\begin{equation*}
\widehat{\mu}(\lambda)=\sum\limits_{g \in G} \mu(g) \lambda(g). 
\end{equation*}
Then its inverse transform is
\begin{equation*}
\mu(g)= \sum\limits_{ \lambda \, \textnormal{irred.}}d_{\lambda} tr(\lambda(g^{-1})\widehat{\mu}(\lambda)),
\end{equation*}
where the sum is over all irreducible representations of $G$ and $d_{\lambda}$ stands for the dimension of the representation $\lambda.$ The inverse transform leads to 
Plancherel's formula stated below.
\begin{equation} \label{planc}
\sum\limits_{g \in G} |\mu(g)|^2 = \frac{1}{|G|} \sum\limits_{ \lambda \, \textnormal{irred.}}d_{\lambda} tr(\widehat{\mu}(\lambda)\widehat{\mu}(\lambda)^{T}).
\end{equation}
It establishes the connection with the total variation distance if we take $\mu$ on the left hand side to be the difference of two measures.

We also note a fact about the Fourier transform of convolutions of measures, as a Markov chain at step $t$ can be viewed as $t-$fold convolution measure of the transition probabilities multiplied by the intial state vector. The convolution of two measures $\mu$ and $\nu$ is defined to be
\begin{equation*}
\mu * \nu (g) = \sum\limits_{h \in G} \mu(g h^{-1}) \nu(h). 
\end{equation*}
The Fourier transform of a convolution satisfies 
\begin{equation} \label{conv}
\widehat{\mu * \nu}(g)= \widehat{\mu}(g) \widehat{\nu}(g). 
\end{equation}

\subsection{Young-Jucys-Murphy elements and eigenvalues} \label{YJMsect} We will focus on the symmetric group $S_n.$ The irreducible representations of $S_n$ are indexed by partitions of integers. We use the notation $\lambda=(\lambda_1,\lambda_2,\ldots) \vdash n$ for partitions of $n.$ $\lambda'$ refers to the conjugate partition of $\lambda$. The dimension of an irreducible representation $\lambda,$ $d_{\lambda},$ is equal to the number of standard Young tableaux of shape $\lambda.$ See, for instance, Chapter 2 of \cite{Sag}. 

First we state an upper bound theorem, which was proved in Chapter 3B of \cite{D} using Plancherel formula \eqref{planc}. As we identify the eigenvalues, we construct the upper bound at the same time. We first note that the irreducible representations of the symmetric group $S_n$ are indexed by partitions of $n.$   and    See Chapter 2 of \cite{Sag}
\begin{lemma} \cite{D}
Let $\mu$ be a probability distribution over $S_n$ and $\pi$ be the uniform distribution. Then,
\begin{equation*}
\| \mu - \pi \|_{TV}^2 \leq \, \frac{1}{4} \sum_{\substack{\lambda \vdash n \\ \lambda \neq (n)}} d_{\lambda} \, tr(\widehat{\mu}(\lambda)\widehat{\mu}(\lambda)^T).
\end{equation*}
\end{lemma}
In our case, we replace $\mu$ by $P_{\theta}^{t}$. The first of the two observations to be made on $P_{\theta}^{t}$, in order to simplify the right hand side, is $P_{\theta}^{t}$ is $t-$fold convolution of $P_{\theta}$. Therefore, by \eqref{conv}, $\widehat{P_{\theta}^{t}}(\lambda)=\widehat{P_{\theta}}(\lambda)^t.$ 

For the second one, first consider $\widehat{P_{\theta}}$ as an element of the group algebra, $\mathbb{C}(S_n).$ A well-known fact about the center of the group algebra, $\mathbf{Z}(\mathbb{C}(S_n))$, is that it has basis consisting of the sum of the all elements in the same conjugacy class. Since the density function defined by $P_{\theta}$ has an equal weight for the permutations in the same conjugacy class, $\widehat{P_{\theta}}(\lambda)$ is in $\mathbf{Z}(\mathbb{C}(S_n))$. But then $\widehat{P_{\theta}}$ is a constant times the identity matrix of the same dimension with $\lambda$, since the representation of any element in the center is in that form as a result of Schur's lemma. We have more to infer from this fact, but for now we use below the implication $\widehat{P_{\theta}}(\lambda)=\widehat{P_{\theta}}(\lambda)^T$.

Putting two observations together, the upper bound can be expressed as
\begin{equation} \label{UB}
\| P_{\theta}^t - U \|_{TV} \leq \, \frac{1}{4} \sum_{\substack{\lambda \vdash n \\ \lambda \neq (n)}} d_{\lambda} tr(\widehat{P_{\theta}}(\lambda)^{2t}).
\end{equation}  
  
Next, we evaluate $\widehat{P_{\theta}}(\lambda).$ We already stated that $\widehat{P_{\theta}}(\lambda)$ is a constant multiple of the identity matrix. In order to identify the constant we first introduce Young-Jucys-Murphy elements. The detailed treatment of the subject can be found in \cite{C,OV}. For $i=2,...,n$, consider the following elements $R_i \in \mathbb{C}(S_n)$:
\begin{equation} \label{YJM}
R_i=(1,i)+(2,i)+\cdots+(i-1,i).
\end{equation}
These elements form a basis for Gelfand-Tsetlin subalgebra of $\mathbb{C}(S_n)$, which is generated by the centers of the symmetric groups $S_1,..,S_n.$ It is a maximal commutative subalgebra, and is the algebra of all operators diagonal in the Gelfand-Tsetlin basis. The diagonal entries of their irreducible representations are identified in \cite{Mu} to be the \textit{contents} of the standard Young tableaux. The content of a box in the Young diagram of $\lambda \vdash n$ is defined by its position in the diagram, it is the column number of the box minus the row number of it. In particular, the upper-left diagonal of $\lambda(R_i)$ is determined by the position of $i$ in the standard Young tableaux which have entries according to the lexicographic order, i.e., $1$ through $\lambda_1$ in the first part, $\lambda+1$ through $\lambda_1+\lambda_2$ in the second part and so on. Denote the upper-left corner element of $\lambda(R_i)$ by $c_{\lambda}(i).$ So the $d_{\lambda} \times d_{\lambda}$ diagonal matrix $\lambda(R_i)$ is in the form

\begin{center}
$\lambda(R_i)=$\(
  \begin{pmatrix}
   \, \, c_{\lambda}(i) & & \\
   & *  & \\
    && \ddots & \\
    & & & *
  \end{pmatrix}
\) 
\end{center}
Although no $R_i$ for $i \in \{2,...,n\}$ is in the center of the group algebra, any symmetric polynomial with variables $R_2,...,R_n$ is in $\mathbf{Z}(\mathbb{C}(S_n)).$ See \cite{G} for the details. We are particularly interested in the elementary symmetric functions of 
YJM elements, the reason is the following theorem:
\begin{theorem}\cite{DG} \label{DGth}
Let $R_i$ be defined as in (\ref{YJM}). The $k^{th}$ elementary symmetric function in $R_2,...,R_n$ is the sum over all permutations in conjugacy classes of $S_n$ with $n-k$ cycles in the group algebra. i.e.,
\begin{equation*}
e_k(R_2,...,R_n) \equiv \sum_{2\leq i_1 <\cdots<i_k\leq n} R_{i_1}\cdots R_{i_k} =  \sum_{\substack{\sigma \in S_n \notag \\ \gamma(\sigma)=n-k}}\sigma.
\end{equation*}
\end{theorem}
A useful corollary, which is related to our observations on the representations of $R_i$, is as follows:
\begin{corollary} \label{DGcor}
Let $\lambda$ be an irreducible representation of $S_n.$ Then
\begin{equation*}
\lambda \left( e_k(R_2,...,R_n) \right) = e_k(c_{\lambda}(2),...,c_{\lambda}(n)) \, I_{d_{\lambda}}.
\end{equation*} 
\end{corollary}
\pf Since $e_k(R_2,...,R_n)$ is the sum of elements in the same conjugacy classes, it is in the center of the group algebra. Therefore, its representation is a constant multiple of the identity matrix by Schur's lemma. So it suffices to consider the upper-left corner entries of each $\lambda(R_i)$ for $i \in \{2,...,n\}$, which are $c_{\lambda}(i)'$s.   \bbox

Now we are ready to calculate $\widehat{P_{\theta}}(\lambda)$ By Theorem \ref{DGth} and its corollary,
\begin{align*}
\widehat{P_{\theta}}(\lambda)&= \sum_{\sigma \in S_n} P_{\theta}(\sigma) \lambda(\sigma) \notag \\
&=\sum_{k=0}^{n-1}\frac{ \theta^{n-k}}{\theta^{(n)}}\sum_{\substack{\sigma \in S_n \notag \\ \gamma(\sigma)=n-k}} \lambda(\sigma) \\
& = \sum_{k=0}^{n-1} \frac{ \theta^{n-k}}{\theta^{(n)}} \, \lambda \left( e_k(R_2,...,R_n) \right) \\
& = \Bigg(\sum_{k=0}^{n-1}\frac{ \theta^{n-k}}{\theta^{(n)}} \, e_k(c_{\lambda}(2),...,c_{\lambda}(n)) \, \Bigg) I_{d_{\lambda}}  \\
& = \Bigg(\frac{ \theta^{n}}{\theta^{(n)}} \sum_{k=0}^{n-1} e_k \left(\frac{c_{\lambda}(2)}{\theta},...,\frac{c_{\lambda}(n))}{\theta}\right) \, \Bigg) I_{d_{\lambda}}.  
\end{align*}
Next, we use the identity
$$ \prod_{i=1}^n (1+x_{i})= \sum_{k=1}^n e_k(x_1,...,x_n) $$
to obtain above,
\begin{align*}
\widehat{P_{\theta}}(\lambda)&=  \Bigg( \frac{ \theta^{n}}{\theta^{(n)}} \, \prod_{i=2}^n \left(1+\frac{c_{\lambda}(i)}{\theta}\right)  \Bigg) I_{d_{\lambda}} \\
&=  \Bigg( \frac{ \theta}{\theta^{(n)}} \, \prod_{i=2}^n (\theta+c_{\lambda}(i))  \Bigg) I_{d_{\lambda}} \\
\end{align*}
Since $c_{\lambda}(1)=0$ for all $\lambda \vdash n$, the expression eventually simplifies to
\begin{equation*}\label{cont}
\widehat{P_{\theta}}(\lambda)=\Bigg(  \prod_{i=1}^{n} \left(\frac{\theta+c_{\lambda}(i)}{\theta+i-1}\right)  \Bigg) I_{d_{\lambda}}. \\
\end{equation*}

Since $P_{\theta}$ is a symmetric matrix, the eigenvalues of irreducible representations give complete list of eigenvalues of $P_{\theta}.$ See  Section 4 of \cite{DS} for the details. Hence, the eigenvalues of $P_{\theta}$ with their multiplicities are
\begin{equation} \label{eigen}
\beta_{\lambda, \theta}=  \prod_{i=1}^{n}\left(\frac{\theta+c_{\lambda}(i)}{\theta+i-1}\right)  \text{ with multiplicity } d_{\lambda}^2 \text{ for } \lambda \vdash n.
\end{equation}

\begin{figure}[H]
\ytableausetup{nobaseline}
\ytableausetup{centertableaux}
\ytableausetup{boxsize=3.5em}
\begin{ytableau}
1& 1 & 1& \ytableausetup{tabloids} \dots  & &  & 1 \ytableausetup{notabloids} &  1\\  \scalebox{1}{$\frac{\theta-1}{\theta+\lambda_1}$} & \scalebox{0.9}{$\frac{\theta}{\theta+\lambda_1+1}$}& \ytableausetup{tabloids} \dots &  &  \ytableausetup{notabloids}  \scalebox{0.8}{$\frac{\theta+\lambda_2-1}{\theta+\lambda_1 + \lambda_2-1}$} \\
\frac{\theta-2}{\theta+\lambda_1 + \lambda_2} & \scalebox{0.8}{$\frac{\theta-1}{\theta+\lambda_1 + \lambda_2 + 1}$} & \none[\dots] \\
\vdots & \none[\vdots] \\
\frac{\theta-\lambda_1'+1}{\theta+n-\lambda_k} & \none[]
\end{ytableau}
\hspace*{1cm}
\caption{Factors of $\beta_{\lambda}$, the eigenvalue associated with $\lambda=(\lambda_1,\ldots,\lambda_k) \vdash n.$}
\end{figure}

\indent Therefore the upper bound on the total variation distance (\ref{UB}) is reduced to
\begin{align}
d_{\theta}(t) & \leq \, \frac{1}{4} \sum_{\substack{\lambda \vdash n \\ \lambda \neq (n)}} d_{\lambda} \, tr\Bigg( \Bigg( \prod_{i=1}^{n} \left(\frac{\theta+c_{\lambda}(i)}{\theta+i-1}\right)^{2t} \Bigg)  I_{d_{\lambda}}\Bigg) \notag \\
& \leq \, \frac{1}{4} \sum_{\substack{\lambda \vdash n \\ \lambda \neq (n)}} d_{\lambda}^2 \,   \prod_{i=1}^{n} \left(\frac{\theta+c_{\lambda}(i)}{\theta+i-1}\right)^{2t}. \label{upperb}
\end{align} 

\section{Proof of Theorem \ref{main} Part (i)}

In this section, we prove the first part of the theorem. It can be explicitly stated that the convergence (to their stationary distributions) times of the Markov chains defined on $S_n$ by the Ewens distributions with a fixed parameter $\theta$ is bounded above by a constant, which depends on only $\theta$, for all $n.$ 

First, let us consider

\begin{equation}\label{l2}
Z_n(\theta):= \sum_{\lambda \vdash n} d_{\lambda}^2 \beta_{\lambda, \theta}^2,
\end{equation}
which is none other than the sum of the eigenvalues of $P_{\theta}^2.$ We will relate this expression to the upper bound on the mixing time. Observe that unlike \eqref{upperb}, the trivial partition $(n)$ is included in the sum. We will calculate $Z_n(\theta)$ explicitly. First we need some definitions.

\subsection{Plancherel growth process}

Let $\mathbb{Y}$ denote the \textit{Young's lattice}, which is the lattice of all integer partitions. A partial order on $\mathbb{Y}$ is defined by the inclusion of Young diagrams, i.e., for $\lambda, \mu \in \mathbb{Y}$ we say $\mu$ is covered by $\lambda,$ and denote it by $\mu \nearrow \lambda,$ if $\lambda$ can be obtained from $\mu$ by adding a box to an admissible position in its Young diagram. This box adding process gives 
a Markov chain on $\mathbb{Y},$ known as the Plancherel growth process, see \cite{K}. Its transition probabilities are given by
\begin{equation*}
p(\lambda, \Lambda) = \frac{d_{\Lambda}}{|\Lambda| d_{\lambda}} \textnormal{ for } \lambda \nearrow \Lambda. 
\end{equation*}
The law of the chain at stage $n$ is the Plancherel measure on $S_n,$ which is defined as
\begin{equation*}
\mathcal{P}(\lambda)=\frac{d_{\lambda}^2}{n!} \textnormal{ for } \lambda \vdash n.
\end{equation*}
We use this process to compute \eqref{l2} recursively. Let us list some properties that the transition measure satisfies. Firstly,

\begin{equation}\label{growth}
d_{\Lambda} = \sum_{\lambda: \lambda \nearrow \Lambda} d_{\lambda} \quad \textnormal{ and } \quad d_{\lambda} = \frac{1}{|\Lambda|} \sum_{\Lambda: \lambda \nearrow \Lambda} d_{\Lambda}.
\end{equation}

Secondly, we state the relation of this measure to contents, similar to the way it is presented in Section 10 of \cite{K2}. Let $X$ be the content of the random box added to the diagram of $\lambda$ according to the transition measure. We have 
\begin{equation}  \label{rancont}
\begin{aligned}
\mathbf{E}_{\lambda}(X)&=\sum_{\Lambda: \lambda \nearrow \Lambda} c(\Lambda \setminus \lambda) \frac{d_{\Lambda}}{|\Lambda|d_{\lambda}}=0, \\
\mathbf{E}_{\lambda}(X^2)&= \sum_{\Lambda: \lambda \nearrow \Lambda} c(\Lambda \setminus \lambda)^2 \frac{d_{\Lambda}}{|\Lambda|d_{\lambda}}=|\lambda|.
\end{aligned}
\end{equation}

\begin{lemma} If $Z_n(\theta)$ be defined as in \eqref{l2}, then 
\begin{equation*}
Z_n(\theta)= \binom{n+\theta^2-1}{n} \bigg/\binom{n+\theta-1}{n}^2.
\end{equation*}
\end{lemma}
\pf By the first equation in \eqref{growth}, we have 
\[Z_n(\theta)=\sum_{\lambda \vdash n} d_{\lambda}^2 \beta_{\lambda,\theta}^{2} = \sum_{\mu \vdash n-1} d_{\mu}^2 \beta_{\mu,\theta}^{2}\sum_{\lambda : \mu \nearrow \lambda} \frac{d_{\lambda}}{ d_{\mu}} \left(\frac{\beta_{\lambda,\theta}}{\beta_{\mu,\theta}}\right)^{2} \]
Next, using the formula \eqref{eigen} and the equations in \eqref{rancont},

\begin{align*}
Z_n(\theta)&= \sum_{\mu \vdash n-1} d_{\mu}^2 \beta_{\mu}^{2}\sum_{\lambda : \mu \nearrow \lambda} \frac{d_{\lambda}}{ d_{\mu}} \left(\frac{\theta+ c(\lambda \setminus \mu)}{n+\theta-1}\right)^{2} \\
& = \frac{n}{(n+\theta-1)^2}\sum_{\mu \vdash n-1} d_{\mu}^2 \beta_{\mu}^{2} \left( \mathbf{E}(\theta^2)+ \mathbf{E}(2 \theta X) + \mathbf{E}(X^2) \right) \\
& = n\frac{n+\theta^2-1}{(n+\theta-1)^2}Z_{n-1}(\theta)
\end{align*}

Since $Z_1(\theta)=1,$ the recursive relation above gives the result.

\bbox

We also note the following combinatorial identity. Since $Z(\theta)$ is in fact the sum of the eigenvalues of $P_{\theta}^2,$ it will agree with the trace of $P_{\theta}^2.$ Each element along the diagonal of $P_{\theta}^2$ is obtained by the sum of squares of elements in any row since $P_{\theta}$ is symmetric. Putting these observations together, we obtain

\[\frac{1}{n!}\sum_{\sigma \in S_n} \theta^{2 \gamma(\sigma)}=\binom{n+\theta^2-1}{n}.\]

If we let $Y(\pi)$ be the number of cycles in a uniformly chosen permutation $\pi \in S_n,$ and define $\Gamma_n(\pi) = \theta^{Y(\pi)},$ we have
\begin{align*}
\mathbf{E}(\Gamma_n)= \binom{n+\theta-1}{n} \quad \textnormal{ and } \quad \mathbf{E}(\Gamma_n^2)& = \binom{n+\theta^2-1}{n}.
\end{align*}
The first identity above is obtained by computing the trace of $P_{\theta}.$

\subsection{Upper bound} Before evaluating the upper bound on the total variation distance, we give an upper on the eigenvalues.
In particular we show that the eigenvalue associated with $\lambda=(n-1,1)$ is the largest among all but the trivial partition $(n)$. 

\begin{lemma}\label{eigmon}
Let $\lambda \neq (n)$ be a partition of $n.$ Then $\beta_{(n-1,1), \theta}^2 \geq \beta_{\lambda, \theta}^2$ for all $\theta >1.$
\end{lemma}

\pf We will compare the factors of the eigenvalues $\beta_{\lambda,\theta}$ and $\beta_{(n-1,1),\theta}$ going through all boxes in their diagrams according to lexicographic order. We ignore the product of denominators in \eqref{eigen}, $\theta^{(n)},$ as its common to all partitions. Every partition has a box in the upper left corner of its diagram, so the first factors agree. Since $\lambda$ is not the trivial partition, its diagram also includes a box in the leftmost position of its second row, which we cross out with the box of $(n-1,1)$ in the same position. Finally, we distinguish two cases, the positive (1) and the negative (2) factors of $\beta_{\lambda,\theta}.$ Since $\theta>1,$ $\beta_{(n-1,1),\theta}$ does not have a negative factor.

Skipping the aforementioned two boxes, we start with the box in the first row and the second column of $(n-1,1).$ It has the factor $(n+\theta+2-1)=\theta+1.$ While the partition $\lambda$ can have a box at the same position or a box at the second row  and the second column or at the third row and the first column. In the last two cases, the factor is $\theta-1,$ which is positive and smaller than of $(n-1,1).$ We repeat this argument for all boxes with positive factors. 

Suppose the $k$th box of $\lambda$ has a negative factor, i.e., $\theta+c_{\lambda}(k+2) < 0.$ Since $c_{\lambda}(k+2)$ is minimized if the first $k$ boxes of $\lambda$ are located in th first column, we have $c_{\lambda}(k+2) > -(k+1).$ On the other hand, the $k$th factor for the partition $(n-1,1)$ is $\theta+k.$  Since $\theta>1,$ we have $|\theta+k| \geq |\theta + c_{\lambda}(k+2)|,$ which proves the second case and the lemma.

\bbox

We are ready to take on \eqref{UB},

\begin{equation}
\begin{aligned}
d_{\theta}(t) & \leq \sum_{\substack{\lambda \vdash n \\ \lambda \neq (n)}}d_{\lambda}^2 \beta_{\lambda}^{2(t+1)} \\
 &\leq \left(\max_{\substack{\lambda \vdash n \\ \lambda \neq (n)}} \beta_{\lambda}^2\right)^{t} \left(Z(\theta)-1 \right) \\
 & \leq \left(\frac{\theta-1}{\theta+n-1}\right)^{2t} \binom{n+\theta^2-1}{n} \binom{n+\theta-1}{n}^{-2}.
\end{aligned}
\end{equation}
Next, we bound the binomial terms in the sum. Let us take $k \leq \theta <k+1$ for $k \in \mathbb{Z}^{+}.$ For the first term, we have
\begin{equation*}
\binom{n+\theta^2-1}{n} < \frac{n^{(k+1)^2}}{k^2!}.
\end{equation*}
We need a lower bound for the latter,
\begin{equation*}
\binom{n+\theta-1}{n}  > \frac{n^{k}}{(k+1)!}.
\end{equation*} 

So, if we take $t=\alpha\frac{k^2+1}{2}$ for some $a>1,$ 

\begin{align*}
d_{\theta}(t)  &\leq  \frac{1}{4} \, \frac{n^{(k+1)^2}}{k^2!} \frac{(k+1)!^2}{n^{2k}} \left( \frac{k -1}{ k+n -1}\right)^{a(k^2+1)} \\
& =C(k)n^{1-\alpha}
\end{align*} 
where $C(k)$ is a some constant depending only on $k=\floor{\theta}.$ Since $\alpha$ is also independent of $n,$  we can take any fixed value for $n,$ say $n\geq 3,$ and conclude that $d_{\theta}(t) < C e^{-\alpha}.$

\bbox

\subsection{A second proof for integer values of the parameter}

Let $k \in \mathbb{Z}^+\setminus\{1\}.$ By hook-length formula, we can write $Z_n(k)$ explicitly as

\begin{equation}
Z_n(k)= \sum_{\lambda \vdash n} \frac{(n!)^2}{(\prod_{i=1}^n h(i))^2} \,   \prod_{i=0}^{n-1} \left(\frac{k+c_{\lambda}(i+1)}{k+i}\right)^{2t},
\end{equation}
Then, we can rearrange it as the following:
\begin{equation}\label{rearr}
\begin{aligned}
Z_n(k)&= \sum_{\lambda \vdash n}  \left(\frac{n!}{\prod_{i=1}^n {k+i-1}}\right)^2 \, \left( \prod_{i=1}^n \frac{k+c_{\lambda}(i)}{h(i)}\right)^2  \\
& = \binom{n+k-1}{n}^{-2} \sum_{\lambda \vdash n}  s_{\lambda}^2(1^k)
\end{aligned}
\end{equation}
where $s_{\lambda}(1^k)= \prod_{i=0}^n \frac{k+c(i)}{h(i)}$ is the Schur function with $k$ non-zero variables, each being equal to $1.$ Note that $s_{\lambda}(1^k)$ is zero unless $\lambda_1' \leq k.$  See Section I-3 of \cite{M} for details. 

Consider the Cauchy identity,
\begin{equation*}
\prod_{i,j} (1-x_i y_j)^{-1} = \sum_{\lambda \text{ part.}} s_{\lambda}(\mathbf{x}) s_{\lambda}(\mathbf{y})
\end{equation*}
where $\mathbf{x}=(x_1,x_2,...),$ $\mathbf{y}=(y_1,y_2,...)$ and the sum is over all partitions $\lambda$ of any size. The proof is found in Section I-4 of \cite{M}. We are interested in partitions $\lambda$ such that $|\lambda|=n,$ so we consider the coefficient of $t^n$ for the generating function 
\begin{equation*}
F(\mathbf{x}, \mathbf{y};t) = \prod_{i,j} (1-x_i y_j t )^{-1}.
\end{equation*}
Then, we take $x_1=\cdots=x_k=1$ and $x_l=0$ for $l > k.$ We do the same for $\mathbf{y}$ to obtain
\begin{equation*}
\sum_{\lambda \text{ part.}} s^2_{\lambda}(1^k) = \, \prod_{i=1}^k \prod_{j=1}^k (1-t)^{-1}= (1-t)^{-k^2}.
\end{equation*}
The right hand side is the generating function for complete symmetric functions, and the coefficient of $t^n$ is given by $h_n=\binom{n+k^2-1}{n}$.
Therefore,
\begin{equation*}
 \sum_{\lambda \vdash n} s^2_{\lambda}(1^k)  = \binom{n+k^2-1}{n}.
\end{equation*}
We then can proceed as before and show the same upper bound. The reason that we cannot simply extend the result from integers to continuum is the failure of eigenvalue monotonicity in the non-integer case. Observe that the eigenvalues are always non-negative for integer values of $\theta,$ which is not true in general. 

\section{Proof of Theorem \ref{main} Part (ii)}\label{pf2}

\subsection{Defining representation} It follows from the definition of the total variation distance (\ref{TV}) that
\begin{equation} \label{lb}
d_{\theta}(t) \geq |P_{\theta}^t(A) - \pi(A)|
\end{equation}
for any subset $A$ of $S_n$. A statistic commonly used to choose an expedient subset in \eqref{lb} is the number of fixed points. Define $A_k$ to be the set of permutations with less than or equal to $k$ fixed points. This particular choice of subsets establishes an important connection with the symmetric group representations.

The defining representation of $S_n$ is the $n$-dimensional representation $\rho$ where 
\begin{equation*}
\rho(\sigma)(i,j)=
 \left\{
\begin{array}{ll}
      1& \sigma(j)=i \\
      0 & \textrm{otherwise} \\
\end{array} 
\right. 
\end{equation*}
for $\sigma \in S_n.$ Denote the character of the representation, or the traces of the matrices given above, by $\chi_{\rho}.$ So that $\chi_{\rho}(\sigma)$ counts the number of fixed points of $\sigma \in S_n.$ Therefore $\sigma \in A_k$ if and only if $\chi_{\rho}(\sigma) \leq k.$

Next, we consider $\chi_{\rho}$ as a random variable under the law of the Markov chain at step $t.$ The defining representation is reducible and decomposed as
\begin{equation} \label{fmom}
\rho=S^{(n)} \oplus S^{(n-1,1)}
\end{equation}
where $S^{\lambda}$ is the Specht module associated with partition $\lambda$ of $n$, noting that the partitions of $n$ are in one to one correspondance with the irreducible representations of $S_n.$ See Chapter 2 of \cite{Sag} for details. For the second moment of $\chi_{\rho},$ we consider the decomposition of $\rho \otimes \rho$ as
\begin{equation}\label{smom}
\rho\otimes \rho =2 S^{(n)} \oplus 3 S^{(n-1,1)} \oplus S^{(n-2,2)} \oplus S^{(n-2,1^2)}.
\end{equation}
Then we use the facts below,
\begin{equation*}
 \begin{split}
\chi_{\rho_1 \oplus \rho_2}&=\chi_{\rho_1}+ \chi_{\rho_2}, \\
\chi_{\rho_1 \otimes \rho_2}&=\chi_{\rho_1} \cdot \chi_{\rho_2},
 \end{split}
\end{equation*}
and invoke some Fourier analytic results, which can be found in \cite{B} Chapter 16 with detailed proofs and in \cite{D} Section 2C in most relevance to our case, to calculate the first two moments of $\chi_{\rho}$. We have the following expressions for the first two moments of a distribution $\mu$ over $S_n$.
\begin{equation}\label{mu}
\begin{aligned}
\mathbf{E}_{\mu}(\chi_{\rho}) &= \text{tr}(\hat{\mu}(n))+ \text{tr}(\hat{\mu}(n-1,1)), \\
\mathbf{E}_{\mu}(\chi_{\rho}^2) &= 2\,\text{tr}(\hat{\mu}(n))+ 3 \, \text{tr}(\hat{\mu}(n-1,1)) +  \text{tr}(\hat{\mu}(n-2,2)) + \text{tr}(\hat{\mu}(n-2,1,1)). 
\end{aligned}
\end{equation}
Then using the first two moments, the probabilities can be bounded by  Chebyshev's inequality,
\begin{equation*}\label{cheby}
\textbf{Pr}(\chi_{\rho} \leq \mathbf{E}(\chi_{\rho}) -a ) \leq \frac{\text{Var}(\chi_{\rho})}{\text{Var}(\chi_{\rho})+a^2}.
\end{equation*}
In particular, the distribution of $\chi_{\rho}$ under the uniform distribution over $S_n$ is well-known; it is given below by the classical matching problem.
\begin{equation} \label{feller}
\pi(A_k)=\mathbf{Pr}_{\pi}(\chi_{\rho} \geq k)=\frac{1}{(k-1)!} \sum_{l=k}^n \frac{(-1)^{l-k}}{l(l-k)!},
 \end{equation}
 for $i=1,2,\dots,$  for which the limiting distribution is Poisson with parameter one. See \cite{T} for various derivations of this result.

\subsection{Lower bound}

The decompositions (\ref{fmom}) and (\ref{smom}) of the defining representation allow us to evaluate the expression in (\ref{mu}) with respect to the law of $P_{\theta}^t.$ We first identify the traces of the Fourier transform by our previous derivation $(\ref{cont})$ for $\widehat{P_{\theta}^t}(\lambda)$ for $\lambda \vdash n.$
\begin{align*}
\text{tr}(\widehat{P_{\theta}^t}(n)) &= 1, \\
\text{tr}(\widehat{P_{\theta}^t}(n-1,1)) &=(n-1)\left(\frac{\theta -1}{\theta + n -1} \right)^t, \\
\text{tr}(\widehat{P_{\theta}^t}(n-2,2))  &= \frac{n (n-3)}{2}\left(\frac{\theta -1}{\theta + n -1} \right)^t \left(\frac{\theta -2}{\theta + n -2} \right)^t, \\
\text{tr}(\widehat{P_{\theta}^t}(n-2,2))  &= \frac{(n-1)(n-2)}{2}\left(\frac{\theta -1}{\theta + n -1} \right)^t \left(\frac{\theta}{\theta + n -2} \right)^t. 
\end{align*}

If we take $t= \frac{\log n}{\log (\theta + n) - \log (\theta)} - \frac{\gamma}{\log (\theta + n) - \log (\theta)}$ where  $\gamma > 0,$ then 
\begin{align*}
\lim_{n \rightarrow \infty}\mathbf{E}_{P_{\theta}^t}(\chi_{\rho})& = 1+e^{\gamma}\\
 \lim_{n \rightarrow \infty}\mathbf{E}_{P_{\theta}^t}(\chi_{\rho}^2)& = 2+3(1+e^{\gamma})+e^{2\gamma}.
\end{align*}
Therefore, $\text{Var}_{P_{\theta}^t}(\chi_{\rho})$ is asymptotically $4+e^{\gamma}. $

We take $n$ large enough and define $B$ to be the set of permutations with less than or equal to $e^{\gamma-1}+1$ fixed points i.e.,

\begin{equation*}
B\equiv\{\sigma \in S_n : \chi_{\rho}(\sigma) \leq e^{\gamma-1}+1 \}.
\end{equation*}
 Then by Chebyshev's inequality \eqref{cheby},
\begin{equation*}
P_{\theta}^t(B)=\textbf{Pr}_{P_{\theta}^t}(\chi_{\rho} \leq 1+e^{\gamma-1}) \leq \frac{4+ e^{\gamma}}{4+ e^{\gamma} + e^{2 \gamma -2}} < \frac{1}{e^{\gamma-2}}.
\end{equation*}
For the uniform distribution, we have 
\begin{equation*}
\pi(B)= \textbf{Pr}_{\pi}(\chi_{\rho} \leq 1+e^{\gamma-1}) \geq 1- \frac{1}{e^{\gamma}!}
\end{equation*}
by  \eqref{feller}. Putting them together,
\begin{align*}
d_{\theta}(t) =&  \max_{S \subseteq S_n}|P_{\theta}^{t}(S) - \pi(S)| \\
\geq&  |P_{\theta}^{t}(B) - \pi(B)| \\
\geq& 1- \frac{1}{e^{\gamma}!}-\frac{1}{e^{\gamma-2}} > 1-\frac{1}{e^{\gamma}}.
\end{align*}
Therefore, we have the lower bound
\begin{equation} \label{lowbd}
d_{n}\left(\frac{\log{n}}{\log 2} -\gamma\right) \geq 1- \frac{1}{e^{\gamma}}.
\end{equation}

\subsection{Dimension bounds and a partial order on partitions}

We present bounds on the sum of $d_{\lambda}^2$'s for certain sets of partitions to be used in the proof of the theorem. A bound for a fixed length of the first row is given by

\begin{lemma}\cite{D} \label{dub} 
Let $|\lambda|=n$ and $\lambda_1$ be the length of the first row of $\lambda.$ Then,
\begin{equation} \label{eqodub}
\sum\limits_{\lambda=(\lambda_1,\dots)} d_{\lambda}^2 \leq \binom{n}{\lambda_1}^2 (n-\lambda_1)!.
\end{equation}
\end{lemma}

The proof in \cite{D} is by bounding the number of standard Young tableaux. First choose $\lambda_1$ elements for the first row, then count all possible standard Young tableaux that can be formed from the remaining cells.
 
Next, by Lemma \ref{dub}, we bound the sum of dimension squares for the partitions with their first row larger than a linear order of $n.$  We use the notation  $f(n)= \mathcal{O}(g(n))$ to mean $\limsup_{n\rightarrow \infty} \left|\frac{f(n)}{g(n)}\right| < \infty,$ and $f(n)= o(g(n))$ for $\limsup_{n\rightarrow \infty} \left|\frac{f(n)}{g(n)}\right|=0.$

\begin{lemma} \label{eqo}
Let $|\lambda|=n$ and $\alpha \in (0,1).$ As $n \rightarrow \infty,$ for all $\epsilon >0,$
\vspace{2mm}
\begin{itemize}
\item[(i)]$\sum\limits_{ \alpha n\leq\lambda_1 } d_{\lambda}^2 = \mathcal{O}(n^{(1-\alpha+ \epsilon)n}),$\\
\item[(ii)]$\sum\limits_{\substack{\alpha n \leq \lambda_1 \\ \frac{\lambda_1}{2} \leq \lambda_2 \leq \lambda_1}} d_{\lambda}^2 = \mathcal{O}(n^{(1-\frac{3}{2}\alpha+ \epsilon)n})$.
\end{itemize}

\end{lemma}
\pf 
The binomial term in \eqref{eqodub} is simply bounded by $4^n$ considering 
\begin{equation*}
\binom{n}{\lambda_1} < \sum_{i=0}^n \binom{n}{k} = 2^n.
\end{equation*}
The second term in \eqref{eqodub} is bounded by the inequality, $(n-\lambda_1)! \leq n^{n-\lambda_1}.$ Combining the two observations, we have
\begin{align*}
\sum\limits_{ \alpha n\leq\lambda_1 } d_{\lambda}^2 &\leq \sum_i \binom{n}{\alpha n +i}^2 (n- \alpha n -i)! \\
&\leq 4^{n} \sum_i (n- \alpha n -i)! \\
& \leq n 4^{n} n^{(1- \alpha)n}  = \mathcal{O}(n^{(1- \alpha + \epsilon)n})
\end{align*}
for all $\epsilon > 0$, which proves part \text{(i)}. \\

To prove part \text{(ii)}, first note that the proof of Lemma \ref{dub} can easily be extended to yield the fact
\begin{equation*}
\sum\limits_{\lambda=(\lambda_1,\lambda_2,\dots)} d_{\lambda}^2 \leq \binom{n}{\lambda_1}^2 \binom{n- \lambda_1}{\lambda_2}^2 (n-\lambda_1-\lambda_2)!.
\end{equation*}
Therefore, similar to part \text{(i)},
 \begin{align*}
\sum\limits_{\substack{\alpha n \leq \lambda_1 \\ \frac{\lambda_1}{2} \leq \lambda_2 \leq \lambda_1}} d_{\lambda}^2 &\leq \sum_{\alpha n \leq \lambda_1} \sum_{\frac{\lambda_1}{2} \leq \lambda_2} \binom{n}{\lambda_1}^2 \binom{n}{\lambda_2}^2 (n-\lambda_1-\lambda_2)! \\
&\leq \sum_{i} \sum_j \binom{n}{\alpha n + i}^2 \binom{n}{\frac{1}{2} \alpha n +j}^2 (n- \frac{3}{2} \alpha n -i - j)! \\   
& \leq n^2 4^{2n} n^{(1- \frac{3}{2}\alpha) n} = \mathcal{O}(n^{(1-\frac{3}{2}\alpha + \epsilon) n}). 
\end{align*}
\bbox

Next, we define a partial order on the partitions of a fixed natural number $n$, known as the \textit{dominance order}, which is to be used in the proof.  
\begin{definition} \cite{M} \label{dom}
Let $\lambda=(\lambda_1, \lambda_2,...)$ and $\mu=(\mu_1, \mu_2,...)$ are partitions of $n.$ We say that $\lambda$ \textit{dominates} $\mu$, and denote by $\lambda \succeq \mu,$ if
\begin{equation*}
\sum\limits_{i=1}^j \lambda_i \geq \sum\limits_{i=1}^j \mu_i \textnormal{  for all j}\geq 1. 
\end{equation*}
If $\lambda \succeq \mu$ and $\lambda \neq \mu,$ then we write $\lambda \succ \mu.$
\end{definition}

\subsection{Upper bound}
We first note that we will drop the parameter in the notation for the eigenvalues. That is to say we will write $\beta_{\lambda}$ instead of $\beta_{\lambda, \theta}$, as we set $\theta$ to be $n$ for the rest of the section. 

The idea is to show that the large eigenvalues have small multiplicities and vice versa. In summary, we first divide $\{\lambda: \lambda \vdash n \}$ into regions with respect to the length of the first two rows, and identify the partitions which dominates every other partition in those regions. Then the bounds on the dimensions given above are employed.

We start by defining two sets of partitions, of which the eigenvalues are computationally managable. The first set consists of partitions where the cells in the Young diagram are stacked up and right as much as possible for a fixed ratio of the length of the first row to the total numbers of boxes. Formally, for a fixed $\alpha \in (0,1),$  
\begin{equation}\label{defa1}
\lambda^{\alpha}_{(1)} \equiv \left(\alpha n,\dots,\alpha n,r n\right) 
\end{equation}
where $\alpha q +r =1 $ for some $q \in \mathbb{N}$ and $ 0 \leq r < \alpha.$ Although the definition makes sense only if $\alpha n$ and $rn$ are integers,  one can choose the closest integers as $n$ tends to infinity and the proof below works properly. 

The set of partitions defined above respects the dominance order (see \ref{dom}) in the sense that
$ \lambda_{(1)}^{\eta} \succeq \lambda_{(1)}^{\zeta}$ if and only if $\eta \geq \zeta.$ In fact, more is true. It is easy to check that 
\begin{equation}\label{dom1}
\lambda_{(1)}^{\alpha} \succeq \lambda \quad \text{ if and only if } \quad \lambda_1 \leq \alpha n. 
\end{equation}

In the second set of partitions, the length of the second row is fixed to be the half of the first row's length if the length of the first row is less than $\frac{n}{2}$, otherwise half of the number of remaining cells. The definition is as follows:
\begin{equation}\label{defa2}
\lambda_{(2)}^{\alpha} \equiv \begin{cases} 
\left(\alpha n, \frac{\alpha}{2}n\dots,\frac{\alpha}{2}n,\tilde{r}n\right) &  \alpha \in \left(0, \frac{1}{2}\right]\\
\left(\alpha n, \frac{1-\alpha}{2}n,\frac{1-\alpha}{2}n \right) & \alpha \in \left(\frac{1}{2},1\right).
\end{cases}
\end{equation}
where $0 \leq \tilde{r} < \frac{\alpha}{2}.$ Similarly, 
\begin{equation}\label{dom2}
\lambda_{(2)}^{\alpha} \succeq \lambda \quad \text{ if } \quad \lambda_1 \leq \alpha n, \, \lambda_2 \leq \min \left\lbrace\frac{\alpha}{2}n, \frac{1-\alpha}{2}n\right\rbrace.
\end{equation}
 
$\{\lambda_{(1)}^{\alpha} \}_{\alpha}$ lies through the upper sides of the triangle in Fig. \ref{fig2}, while $\{\lambda_{(2)}^{\alpha} \}_{\alpha}$ lies strictly below it through the inscribed triangle in the figure. Then the dominance order can be interpreted as $\lambda_{(i)}^{\alpha}$ dominates every partition to the right under the triangle it lies on. 

Next, we estimate the eigenvalues corresponding to those partitions defined above.  
\begin{lemma}\label{lemb1}
Let $|\lambda|=n$ and $\alpha q +r =1$ for some fixed $\alpha$ as in  definition \eqref{defa1}. Then we have, 
\begin{equation*}
\beta_{\lambda_{(1)}^{\alpha}} \leq C \frac{\left(1 + \alpha \right)^{q(1+\alpha)n} (1+r)^{(1+r)n}}{2^{2n}},
\end{equation*}
where $C$ is a constant depending on $\alpha.$
\end{lemma}

\pf We use the formula for the eigenvalues \eqref{eigen} to have
\begin{align*}
\beta_{\lambda_{(1)}^{\alpha}} &= \prod_{i=1}^q \prod_{j=1}^{\alpha n} \frac{n-i+j}{n+(i-1)\alpha n+j-1} \times \prod_{j=1}^{rn} \frac{n-q+j}{2n-rn+j-1}  \\
&\leq \prod_{i=1}^q \prod_{j=1}^{\alpha n} \frac{n+j}{n+(i-1)\alpha n+j} \times \prod_{j=1}^{rn} \frac{n+j}{2n-rn+j} \\
&\leq \prod_{i=1}^q \frac{[(1+\alpha)n]!}{n!} \frac{[(1+(i-1)\alpha)n]!}{[(1+i\alpha)n]!} \times \frac{[(1+r)n]!}{n!} \frac{[(1+q\alpha)n]!}{(2n)!} \\
&= \frac{([(1+\alpha)n]!)^q \, [(1+r)n]!}{(n!)^{q} \, (2n)!}
\end{align*}
For $n$ large enough, applying the formula (see  6.1.38 in \cite{Ab}) 
\begin{equation}\label{stir}
x!=\sqrt{2 \pi} x^{x+\frac{1}{2}} e^{-x + \frac{\xi}{12x}}, \quad \text{ for some } 0<\xi<1,
\end{equation}
to the factorials above, we obtain
\begin{align*}
\beta_{\lambda_{(1)}^{\alpha}}
&\leq e^{\frac{q}{12n}} \sqrt{\frac{(1+\alpha)^q (1+r)}{2}} \frac{\left(1 + \alpha \right)^{q(1+\alpha)n} (1+r)^{(1+r)n}}{2^{2n}} \\
& \leq 2^{2/\alpha} \frac{\left(1 + \alpha \right)^{q(1+\alpha)n} (1+r)^{(1+r)n}}{2^{2n}}.
\end{align*}
\bbox 

In fact, for our purpose, only the cases where either $r=0$ or $r= \frac{\alpha}{2}$ are relevant. For the second set of partitions, we have the following lemma.
\begin{lemma}\label{lemb2}
Let $|\lambda|=n$ and  $\tilde{r}=0$ in the definition \eqref{defa2} for  $\alpha \in \left(0, \frac{1}{2} \right]$. Then there exists $\tilde{q} \in \mathbb{N}$ such that  $\alpha\tilde{q} =2.$  The eigenvalues are bounded as
\begin{equation*}
\beta_{\lambda_{(2)}^{\alpha}}\leq C \left(1+\frac{\alpha}{2} \right)^{\alpha n /2} \beta_{\lambda_{(1)}^{\alpha /2 }},
\end{equation*}
where $C$ is a constant depending on $\alpha.$
\end{lemma}
\pf The eigenvalue formula \eqref{eigen} gives \\
\begin{align*}
\beta_{\lambda_{(2)}^{\alpha}} =& \prod_{i=1}^{\tilde{q}-2} \prod_{j=1}^{\alpha n/2} \frac{n-i+j}{n+(i+1)\frac{\alpha n}{2}+j-1}   \\
=& \prod_{j=1}^{\alpha n/2} \frac{n+ \frac{\alpha n}{2}+j -1 }{n -(\tilde{q}-1)+j} \, \prod_{i=1}^{\tilde{q}-1} \prod_{j=1}^{\alpha n/2} \frac{n-i+j}{n+i\frac{\alpha n}{2}+j-1}  \\
\leq &  \left( \frac{n+ \frac{\alpha n}{2} -1 }{n -(\tilde{q}-1)}\right)^{\alpha n/2} \beta_{\lambda_{(1)}^{\alpha /2 }} \\
= & C \left(1+\frac{\alpha}{2} \right)^{\alpha n /2}  \beta_{\lambda_{(1)}^{\alpha /2 }}
\end{align*}
for some constant $C.$

\bbox

Having established bounds on set of eigenvalues, we then bound the total variation distance for $\theta=n$ using the upper bound formula \eqref{upperb}. We separate the range of the sum into four regions as depicted in Fig. \ref{fig2}. In set notation, they are 
\begin{align*}
R_1=&\{\lambda \vdash n: \lambda_1 \leq 13\}, \\
R_2=&\left\lbrace\lambda \vdash n: 13 < \lambda_1 \leq 3 , \, \, \min \left\lbrace\frac{\lambda_1}{2}, \frac{1- \lambda_1}{2} \right\rbrace \leq \lambda_2\right\rbrace, \\
R_3=&\left\lbrace\lambda \vdash n: 13 < \lambda_1 \leq 3, \, \,  \lambda_2 < \min \left\lbrace\frac{\lambda_1}{2}, \frac{1- \lambda_1}{2} \right\rbrace  \right\rbrace, \\
R_4=&\{\lambda \vdash n: \frac{n}{3} < \lambda_1 < n\}.
\end{align*}

We take $t^*=\frac{\log n}{\log 2}$ throughout the proof. Also note that
\begin{equation}\label{logn}
 \beta_{\lambda}^{2t^*}= \beta_{\lambda}^{\frac{2}{\log 2}\log n}=n^{\frac{2}{\log 2} \log \beta_{\lambda}}.
 \end{equation}

\begin{figure}
\begin{center}
\begin{tikzpicture}[scale=0.75]

\draw[->] (0,0) -- (9,0) node[anchor=north] {$\lambda_1$};
\draw	(0,0) node[anchor=north] {0}
        (1.5,0) node[anchor=north] {$n/13$}
		(4,0) node[anchor=north] {$n/2$}
		(8,0) node[anchor=north] {$n$}
		(5.6,0)node[anchor=north] {$n/3$};
	
\draw	(1,0.4) node{{\Large $R_1$}}
		(4,2.5) node{{\Large $R_2$}}
		(4,0.9) node{{\Large $R_3$}}
         (6.4,0.7) node{{\Large $R_4$}};
\draw[->] (0,0) -- (0,5) node[anchor=east] {$\lambda_2$};
\draw (0,2) node[anchor=east] {$n/4$}
       (0,4) node[anchor=east] {$n/2$};
\draw[thick] (1.5,0) -- (1.5, 1.5);
\draw[dotted] (4,0) -- (4,4);
\draw[thick] (5.6,0) -- (5.6,2.4);
\draw[thick] (1.5,0.75) -- (4,2);
\draw[thick] (4,2) -- (5.6,1.2);
\draw[dotted] (0,0) -- (1.5,0.75);
\draw[dotted] (5.6,1.2) -- (8,0);
\draw[dotted] (0,2) -- (4,2);
\draw[dotted] (0,4) -- (4,4);
\draw[thick] (0,0) -- (4,4);

\draw[thick] (4,4) -- (8,0);

\end{tikzpicture}
\end{center}
\caption{Planar representation of partitions with respect to the length of the first two rows of Young diagram}\label{fig2}
\end{figure}

\subsubsection{Region 1 } 
We show that the eigenvalues are small enough after a linear order of $n,$ so that when multiplied with $d_{\lambda}^2$ the sum is of smaller order than a constant. 

We take $r=0$ in the statement of Lemma \ref{lemb1}, therefore $\alpha=\frac{1}{q}$ for some $q \in \mathbb{N}$. By Lemma \ref{lemb1},
 \begin{align}
\log \beta_{\lambda_{(1)}^{\alpha}}=& \left(q \left(1+ \frac{1}{q} \right) \log \left(1 + \frac{1}{q} \right) - 2 \log 2 \right)n \notag \\
=& \left( \left(1+q\right) \left( \frac{1}{q} - \frac{1}{2q^2} + \frac{1}{3q^3} - \cdots \right) - 2 \log 2 \right) n \notag \\
\leq & \left( 1+ \frac{1}{2q}  -2 \log 2 \right)n \label{oneoverq}
\end{align}
Then we bound \eqref{logn}.
\begin{align*}
\log_n \left(\beta_{\lambda_{(1)}^{\alpha}}\right)^{2t^{*}}=& \frac{2}{\log 2} \log \beta_{\lambda_{(1)}^{\alpha}} \\
= & \left(\frac{2}{\log 2} + \frac{1}{q\log 2 } - 4\right) n \\
\approx & \left(\frac{1.44}{q} -1.11 \right)n
\end{align*}

So if we choose $q$ large enough, or equivalently $\alpha$ small enough, $\left(\beta_{\lambda_{(1)}^{\alpha}}\right)^{2t^{*}} \leq n^{-n}.$ It turns out that the smallest integer $q$ that satisfies the inequality is $13.$ We have,

\begin{align*}
\sum_{\lambda \in R_1} d_{\lambda}^2 \beta_{\lambda}^{2t^*} \leq & \left(\beta_{\lambda_{(1)}^{1/13}}\right)^{2t^*} \sum\limits_{\substack{\lambda \vdash n \\  \lambda_1 \leq n/13}} d_{\lambda}^2 \\ \leq& n^{-1.01n} \sum_{\lambda \vdash n} d_{\lambda}^2 \\ \leq&  n^{-1.01n} n! \\ =& o(n^{-0.01 n}).
\end{align*}

\subsubsection{Region 2 }
First we briefly justify the need for the restriction on the length of the second row of Young diagrams. Consider the partition $\lambda=(\frac{n}{2}, \frac{n}{2}),$ and the corresponding eigenvalue 
\begin{equation*}
\beta_{\left(\frac{n}{2}, \frac{n}{2} \right)}= \prod_{i=1}^{n/2} \frac{n+i-1}{n+ \frac{n}{2}-1} = \frac{\binom{\frac{3}{2}n-1}{\frac{n}{2}}}{\binom{2n-1}{\frac{n}{2}}} = C\frac{\left(\frac{3}{2}\right)^{3n/2}\left(\frac{3}{2}\right)^{3n/2}}{2^{2n}} =C \frac{3^{3n}}{2^{5n}}
\end{equation*}
for some constant $C.$ Then, 
\begin{equation*}
\beta_{\left(\frac{n}{2}, \frac{n}{2} \right)}^{2t^*}= C n^{\frac{2}{\log 2} n[3 \log 3 - 5 \log 2]} \approx C n^{-0.49 n}.
\end{equation*}
On the other hand, $\sum\limits_{\lambda=(\frac{n}{2}, \dots)} d_{\lambda}^2$ is of order at least $\frac{n}{2}!.$ To see this, first fill out the first row by $1,2,\dots,\frac{n}{2},$ then there are exactly $\frac{n}{2}!$ ways to obtain a standard Young diagram, which gives a lower bound. But since $n^{(1/2- \epsilon) n} = o(\frac{n}{2}!),$ $\sum\limits_{\lambda=(\frac{n}{2}, \dots)} d_{\lambda}^2 \, \beta_{\left(\frac{n}{2}, \frac{n}{2} \right)}^{2t^*}$  is of order larger than a constant.  

We restrict our attention to a smaller region, and use part(ii) of Lemma \ref{eqo} to bound the dimensions. First, we choose a sequence of partitions $\{\lambda_{(1)}^{\alpha_i}\}_i$ in this region, for which $\alpha_0=\frac{2}{3},\alpha_1=\frac{1}{2}, \alpha_2=\frac{2}{5}$ and $\alpha_i=\frac{1}{i}$ for $i \geq 3.$ By the dominance order relations \eqref{dom1}, we have
\begin{equation*}
\sum_{\lambda \in R_2} d_{\lambda}^2 \left(\beta_{\lambda}\right)^{2t^*} \leq \sum_i \left(\beta_{\lambda_{(1)}^{\alpha_i}}\right)^{2 t^*}\sum\limits_{\substack{\lambda \in R_2 \\  \alpha_{i+1}n \leq \lambda_1 \leq \alpha_i n}} d_{\lambda}^2. 
\end{equation*}
Then by part(ii) of Lemma \ref{eqo},
\begin{equation*} 
\sum_{\lambda \in R_2} d_{\lambda}^2 \left(\beta_{\lambda}\right)^{2t^*} \leq \sum_i \left(\beta_{\lambda_{(1)}^{\alpha_i}}\right)^{2 t^*} n^{(1-\frac{3}{2}\alpha_{i+1})n}.
\end{equation*}

Next, we bound the eigenvalues. For $i=0,1 \text{ and } 2,$ we have the following calculations by Lemma \ref{lemb1}. 
\begin{align*}
\log \beta_{(1)}^{\alpha_0} &=\log \beta_{\left(\frac{2n}{3}, \frac{n}{3} \right)} \leq \left(\frac{5}{3} \log \frac{5}{3}  + \frac{4}{3} \log  \frac{4}{3}  -2 \log 2 \right)n \leq -0.15 n, \\
\log \beta_{(1)}^{\alpha_1} &=\log \beta_{\left(\frac{n}{2}, \frac{n}{2} \right)} \leq \left(3 \log \frac{3}{2} -2 \log 2 \right)n \leq -0.17 n, \\
\log \beta_{(1)}^{\alpha_2} &=\log \beta_{\left(\frac{2n}{5}, \frac{2n}{5}, \frac{n}{5} \right)} \leq \left( \frac{14}{5}  \log \frac{7}{5}  +  \frac{6}{5} \log  \frac{6}{5}  -2 \log 2 \right)n \leq -0.22 n.
\end{align*}
Further computations using \eqref{logn} yield
\begin{equation} \label{case21}
\sum_{i=0}^2 \left(\beta_{\lambda_{(1)}^{\alpha_i}}\right)^{2 t^*} n^{(1-\frac{3}{2}\alpha_{i+1})n} = o(n^{-0.09n}).
\end{equation}

For $i \geq 3,$ since $\alpha_i=\frac{1}{i},$ we can use the bound \eqref{oneoverq} to obtain  
\begin{equation*}
\log \beta_{(1)}^{\alpha_i} \leq  \left( 1+\frac{1}{2i} - 2 \log 2  \right)n. 
\end{equation*}
Then by \eqref{logn},
\begin{align*}
&\sum_{i \geq 3} \left(\beta_{\lambda_{(1)}^{\alpha_i}}\right)^{2 t^*} n^{(1-\frac{3}{2}\alpha_{i+1})n} \\
 =& \sum_{i \geq 3} \exp \left\lbrace \left( \frac{2}{\log 2} -4 + \frac{1}{i \log 2} + 1 - \frac{3}{2(i+1)}\right)n\log n \right\rbrace .
\end{align*}
It is easy to see that the right hand side is a decreasing function of $i$ for $i\geq 3.$ Therefore, if we plug in $i=3$ and carry out the calculations,
\begin{align*}
\sum_{i \geq 3} \left(\beta_{\lambda_{(1)}^{\alpha_i}}\right)^{2 t^*} n^{(1-\frac{3}{2}\alpha_{i+1})n} &\leq  \sum_{i \geq 3} \exp \left\lbrace \left( \frac{2}{\log 2} -4 + \frac{1}{3 \log 2} + 1 - \frac{3}{8}\right)n\log n \right\rbrace \\
&=o(n^{-0.008n}).
\end{align*}
Combining with \eqref{case21}, we eventually have
\begin{equation*}
\sum_{\lambda \in R_2} d_{\lambda}^2 \left(\beta_{\lambda}\right)^{2t^*} = o(n^{-0.008n}).
\end{equation*}

\subsubsection{Region 3}

This region is treated very similar to Region 2. We consider the second set of partitions \eqref{defa2} this time with the same choice of $\alpha_i's$ as in the previous case. \\

By the dominance order relations \eqref{dom2} mentioned above, we have
\begin{equation*}
\sum_{\lambda \in R_3} d_{\lambda}^2 \left(\beta_{\lambda}\right)^{2t^*} \leq \sum_i \left(\beta_{\lambda_{(2)}^{\alpha_i}}\right)^{2 t^*}\sum\limits_{\substack{\lambda \in R_3 \\  \alpha_{i+1}n \leq \lambda_1 \leq \alpha_i n}} d_{\lambda}^2. 
\end{equation*}
By part(i) of Lemma \ref{eqo},
\begin{equation*}
\sum_{\lambda \in R_3} d_{\lambda}^2 \left(\beta_{\lambda}\right)^{2t^*} \leq \sum_i \left(\beta_{\lambda_{(2)}^{\alpha_i}}\right)^{2 t^*} n^{(1-\alpha_{i+1})n}.
\end{equation*}

For $i=0,$ we can bound the eigenvalue by direct calculations,
\begin{align*}
 \beta_{(2)}^{\alpha_0} =\beta_{\left(\frac{2n}{3}, \frac{n}{6}, \frac{n}{6} \right)} &= \prod_{i=1}^2 \prod_{j=1}^{n/6} \frac{n-i+j}{\frac{5n}{3}+(i-1)\frac{n}{6}+j-1} \\
 &\leq C \frac{\left(\frac{7}{6}\right)^{14n/6} \left(\frac{5}{3}\right)^{5n/3}}{2^{2n}},
\end{align*}
 which follows from the proof of Lemma \ref{lemb1}. Similarly,
 \begin{align*}
 \beta_{(2)}^{\alpha_1} =&\beta_{\left(\frac{n}{2}, \frac{n}{4}, \frac{n}{4} \right)} \leq C \frac{\left(\frac{5}{4}\right)^{5n/2} \left(\frac{3}{2}\right)^{3n/2}}{2^{2n}}, \\
 \beta_{(2)}^{\alpha_2} =&\beta_{\left(\frac{2n}{5}, \frac{n}{5}, \frac{n}{5}, \frac{n}{5} \right)} \leq C \frac{\left(\frac{6}{5}\right)^{18n/5} \left(\frac{7}{5}\right)^{7n/5}}{2^{2n}}.
 \end{align*}
Therefore,
\begin{align*}
\log \beta_{(2)}^{\alpha_0} &\leq \left( \frac{14}{6} \log \frac{7}{6} + \frac{5}{3} \log \frac{5}{2} -2 \log 2 \right)n \leq -0.175 n, \\
\log \beta_{(2)}^{\alpha_1} & \leq \left( \frac{5}{2} \log \frac{5}{4} + \frac{3}{2} \log \frac{3}{2} -2 \log 2 \right)n  \leq -0.22 n, \\
\log \beta_{(2)}^{\alpha_2} & \leq \left( \frac{18}{5} \log \frac{6}{5} + \frac{7}{5} \log \frac{7}{5} -2 \log 2 \right)n  \leq -0.25 n.
\end{align*}
Then we have
\begin{equation} \label{case31}
\sum_{i=0}^2 \left(\beta_{\lambda_{(2)}^{\alpha_i}}\right)^{2 t^*} n^{(1-\alpha_{i+1})n}= o(n^{-0.005n}).
\end{equation}

For $i \geq3,$ we apply Lemma \ref{lemb1} and \ref{lemb2} to obtain
\begin{align*}
\log \beta_{(2)}^{\alpha_i} \leq&  \left( \frac{1}{2i} \log \left(1+\frac{1}{2i}\right) +2i\left(1+\frac{1}{2i}\right) \log \left( 1+\frac{1}{2i} \right) - 2 \log 2  \right)n \\
=& \left( \left(2i + 1 +\frac{1}{2i} \right) \left( \frac{1}{2i} - \frac{1}{8i^2} + \frac{1}{24 i^3}- \frac{1}{64 i^4} + \cdots \right) -2 \log 2 \right)n \\
\leq & \left(1+ \frac{1}{4i} + \frac{5}{24 i^2} -2 \log 2  \right)n
\end{align*}
Then by \eqref{logn},
\begin{align*}
&\sum_{i=3}^{12}  \left(\beta_{\lambda_{(2)}^{\alpha_i}}\right)^{2 t^*} n^{(1-\alpha_{i+1})n} \\
=& \sum_{i =3}^{12} \exp \left\lbrace \left( \frac{2}{\log 2} -4 + \frac{1}{2i \log 2} + \frac{5}{12 i^2 \log 2}+ 1 - \frac{1}{(i+1)}\right)n\log n \right\rbrace .
\end{align*}
Now we observe that
\begin{equation*}
\frac{2}{\log 2}-3 + \frac{5}{12 i^2 \log 2} \leq \frac{2}{\log 2}-3 + \frac{5}{108 \log 2} < 0,
\end{equation*}
for $i \geq 3.$ One can also verify that $ \frac{1}{2i \log 2}  - \frac{1}{(i+1)}$ is an increasing function of $i$ for $i\geq 3.$ Therefore, we can plug in $i=12$ to obtain an upper bound.
\begin{equation*}
\sum_{i =3}^{12}  \left(\beta_{\lambda_{(2)}^{\alpha_i}}\right)^{2 t^*} n^{(1-\alpha_{i+1})n} 
\leq \sum_{i =3}^{12} \exp \left\lbrace \left(  \frac{1}{144\log 2}  - \frac{1}{13}\right)n\log n \right\rbrace = o(n^{-0.06n}).
\end{equation*}
Putting together with \eqref{case31}, we have
\begin{equation*}
\sum_{\lambda \in R_3} d_{\lambda}^2 \left(\beta_{\lambda}\right)^{2t^*} = o(n^{-0.005n}).
\end{equation*}
\subsubsection{Region 4}
The cases covered above allows us to conclude that if $t \geq \frac{2}{\log 2} \log n,$ Region 4 determines the convergence rate. Take $t= \frac{\log n}{\log 2} + \gamma$ for $\gamma > 0.$ \\ 
Let $m=n- \lambda_1.$ Then we have\
\begin{align*}
\beta_{(\lambda_1, n-\lambda_1)} =& \prod_{j=1}^m \frac{n+j-1}{2n-m+j-1} \\
\leq& \left( \frac{n+m}{2n} \right)^m
\end{align*}
Using the fact that $\log(1+x) < x,$
\begin{equation}\label{logest}
\log \beta_{(n-m,m)} \leq m \left[ -\log 2 + \log \left(1+ \frac{m}{n} \right) \right] 
\leq -m \log 2 + \frac{m^2}{n},
\end{equation}
which implies
\begin{align*}
\beta_{(n-m,m)}^{2t} =\beta_{(n-m,m)}^{\frac{2 \log n}{\log 2}} \, \beta_{(n-m,m)}^{2\gamma}& \leq n^{\frac{2m^2}{n \log 2}-2m} \,  \beta_{(n-1,1)}^{2\gamma} \\
& \leq  n^{\frac{2m^2}{n \log 2}-2m} \left(\frac{1}{2} \right)^{2 \gamma }. 
\end{align*}
Therefore, by Lemma \ref{dub} and \eqref{logest}, 
\begin{align*}
\sum_{\lambda \in R_4} d_{\lambda}^2 \, \beta_{\lambda}^{2t}=&  \sum\limits_{2n/3 \leq \lambda_1 \leq n-1} d_{\lambda}^2 \, \beta_{\lambda}^{2t} \\
\leq & \left(\frac{1}{2} \right)^{2 \gamma } \sum_{m=1}^{n/3} \binom{n}{m}^2 m! \,  n^{\frac{2m^2}{n \log 2}-2m} . 
\end{align*}
To bound the factorial term, we employ the fact $\binom{n}{m} \leq \frac{n^m}{m!}$, 
\begin{align*}
\sum_{\lambda \in R_4} d_{\lambda}^2 \, \beta_{\lambda}^{2t} \leq & \,
\frac{1}{4^{\gamma}} \sum_{m=1}^{n/3} \frac{n^{2m}}{m!} n^{\frac{2m^2}{n \log 2 } -2m} \\
=& \, \frac{1}{4^{\gamma}} \sum_{m=1}^{n/3} \frac{n^{\frac{2m^2}{n \log 2 }}}{m!}.
\end{align*}
An application of Stirling's formula gives
\begin{align*}
\sum_{\lambda \in R_4} d_{\lambda}^2 \, \beta_{\lambda}^{2t} \leq &
\, \frac{1}{4^{\gamma}} \left( 1+ \frac{1}{2} +\frac{e}{\sqrt{2 \pi}}\sum_{m=3}^{n/3} \frac{1}{m^{3/2}} \frac{e^{m-1}n^{\frac{2m^2}{n \log 2 }}}{m^{m-1}} + \mathcal{O}(n^{1/n})\right).
\end{align*}
Let
\begin{equation*}
 f(m)=\log \left(\frac{e^{m-1} n^{\frac{2 \, m^2}{n \log 2}}}{m^{m-1}}\right)=(m-1)+\frac{2 m^2 \log n}{n \log 2}  - (m-1) \log m.
\end{equation*}
Next, we show that $f(m)$ is less than $0$ in the range of the sum. Taking derivatives, we have
\begin{equation*}
f''(m)= \frac{4 \log n}{n \log 2} - \frac{1}{m} -\frac{1}{m^2} >0 
\end{equation*}
for $m \geq 3.$ Therefore $f$ is a convex function, so 
\begin{equation*}
f(m) \leq \max \{f(3), f(n/3) \} 
\end{equation*}
 for all  $m \in [3, n/3].$ It is easy to check that both $f(3)$ and $f(n/3)$ are less than zero, hence $f$ is strictly smaller than zero in the range of the sum. As a result of this,

\begin{align*}
\sum_{\lambda \in R_4} d_{\lambda}^2 \, \beta_{\lambda}^{2t} \leq &
\, \frac{1}{4^{\gamma}} \left( 1+ \frac{1}{2} +\frac{e}{\sqrt{2 \pi}}\sum_{m=3}^{n/3} \frac{1}{m^{3/2}} + \mathcal{O}(n^{1/n})\right) \\
& \leq \, \frac{1}{4^{\gamma}} \left( 1+ \frac{1}{2} +\frac{e}{\sqrt{2 \pi}}(\zeta(3/2)-1)\right) \\
& \leq \, \frac{1}{4^{\gamma-1}}. 
\end{align*}
Finally, by \eqref{upperb} we have
\begin{equation*}
d_{n}\left(\frac{\log{n}}{\log 2} +\gamma\right)  \leq \frac{1}{4^{\gamma}}.  
\end{equation*}
Combining with \eqref{lowbd}, we conclude that the chain has a total variation cutoff.

\bbox

\section{Conclusion}

The Plancherel growth process does not seem to be fully exploited in this context. Further work could require estimates on the moments of contents with respect to the transition measure. An explicit formula for the moments is provided in \cite{L04}. The computations in Section 4.2 suggests the following:
\begin{conjecture}  Let $P_{\theta}$ be the Markov chain on $S_n$ defined in \eqref{MC}. 
\begin{itemize}
\item[(i)] If $\theta(n) = \mathcal{O}(n^\beta)$ for some $\beta < 1,$ then the mixing time is bounded above by a constant independent of $n.$ 
 
\vspace*{1mm}

\item[(ii)] If $n = \mathcal{O}(\theta(n))$, then the chain has a total variation cutoff at $t_n=\frac{\log n}{\log(n+\theta(n))-\log \theta(n)}$ with a window of size $\frac{1}{\log(n+\theta(n))-\log \theta(n)}.$
\end{itemize}
\end{conjecture}

To compare it to a well-known example, we expect that $P_{\theta}$ has the same rate of convergence with the random transposition walk on the symmetric group \cite{DS} provided that $\lim_{n\rightarrow \infty} \frac{\theta(n)}{n^2/2}=1.$ 

\section*{Acknowledgement}
The author would like to thank Jason Fulman for the suggestion of the problem.

\bibliographystyle{alpha}
\bibliography{ewens}

\end{document}